\documentclass[12pt]{article}

\usepackage{amsmath, epsfig, cite, lineno}
\usepackage{amssymb}
\usepackage{amsfonts, color}
\usepackage{latexsym}

\numberwithin{equation}{section}

\begin{document}

\begin{center}
{\Large\bf A combinatorial proof of Guo's multi-generalization\\[5pt] of Munarini's identity}
\end{center}

\vskip 2mm \centerline{Dan-Mei Yang}
\begin{center}
{\footnotesize Department of Mathematics, East China Normal University,\\ Shanghai 200062,
 People's Republic of China\\
{\tt dan2004@126.com} }
\end{center}

\vskip 0.7cm \noindent{\bf Abstract.}  We give a
combinatorial proof of Guo's multi-generalization of Munarini's
identity, answering a question of Guo.


\section{Introduction}
\noindent
 Simons \cite{Simons} proved a binomial coefficient identity
which is equivalent to
\begin{equation}\label{eq1}
\sum_{k=0}^n{n \choose k}{n+k \choose
k}(-1)^{n-k}(1+x)^k=\sum_{k=0}^n {n \choose k}{n+k \choose k}x^k.
\end{equation}
Several different proofs of \eqref{eq1} were given by \cite{Chapman,Prodinger,Sun}. Using Cauchy's integral formula as
in \cite{Prodinger}, Munarini \cite{Munarini} obtained the following
generalization
\begin{equation}\label{eq2}
\sum_{k=0}^n{\beta-\alpha+n \choose n-k}{\beta+k \choose
k}(-1)^{n-k}(x+y)^ky^{n-k}=\sum_{k=0}^n {\alpha \choose n-k}{\beta+k \choose k}x^ky^{n-k},
\end{equation}
\noindent where $\alpha,\beta,x$~and $y$~are indeterminates. It is clear that
\eqref{eq2} reduces to \eqref{eq1} when $\alpha=\beta=n$ and
$y=1$.~Shattuck \cite{Shattuck} and Chen and Pang \cite{Chen} provided
two interesting combinatorial proofs of \eqref{eq2}.

Recently, Guo \cite{Guo} obtained the following multinomial
coefficient generalization of \eqref{eq2} as follows:
\begin{equation}\label{eq3}
\sum_{{\bf k}={\bf 0}}^{{\bf n}}(-1)^{|{\bf n}|-|{\bf k}|}{\beta-\alpha+|{\bf n}|
\choose {\bf n}-{\bf k}}{\beta+|{\bf k}| \choose
{\bf k}}({\bf x}+{\bf y})^{{\bf k}}{\bf y}^{{\bf n}-{\bf k}}
=\sum_{{\bf k}={\bf 0}}^{{\bf n}}{\alpha \choose
{\bf n}-{\bf k}}{\beta+|{\bf k}| \choose
{\bf k}}{\bf x^k}{\bf y}^{{\bf n}-{\bf k}},
\end{equation}
where ${\bf n}=(n_1,\ldots,n_m)\in\mathbb{N}^m$, $|{\bf n}|=n_1+\cdots+n_m,$ ${\bf x},{\bf y}\in\mathbb{C}^m$,
the \emph{multinomial coefficient} ${x \choose {\bf n}}$ is defined by
\begin{equation*}{x \choose {\bf n}}=
\begin{cases}
\displaystyle\frac{x(x-1)\cdots(x-|{\bf n}|+1)}{n_1!\cdots n_m!}, & \text{if}~
{\bf n}\in\mathbb{N}^m,\\0, & \text{otherwise,}
\end{cases}
\end{equation*}
and
${\bf b^a}=b_1^{a_1}\cdots b_m^{a_m}$ for ${\bf a}=(a_1,\ldots,a_m)\in\mathbb{N}^m$ and ${\bf b}=(b_1,\ldots,b_m)\in\mathbb{C}^m$.

In this paper we shall give an involutive proof of
\eqref{eq3}, answering a question of Guo \cite{Guo}. Our proof is motivated by Shattuck's work \cite{Shattuck}.

\section{The Involutive Proof}
\noindent Note that both sides of (\ref{eq3}) are polynomials in
$\alpha, \beta,x_1,\ldots,x_m$ and $y_1,\ldots,y_m$. We may only
consider the case of positive integers with $\beta\geq\alpha$. We
first understand the unsigned quantity in the sum of the left-hand
side of (\ref{eq3}). Let $\Gamma=\{a,b_1,\ldots,b_m\}$ be an alphabet.
We construct the weighted words $w=w_1\cdots w_{\beta+|{\bf n}|}$ on $\Gamma$ as follows:
\newcounter{xz}
\begin{list}{\roman{xz})\hfill}%\bfseries\sffamily
            {\setlength{\topsep}{2mm}
             \setlength{\labelwidth}{1em}
             \setlength{\leftmargin}{2em}
             \setlength{\listparindent}{0pt}
             \setlength{\itemsep}{2pt}
             \setlength{\parsep}{5pt}
             \usecounter{xz}}
\item
Choose ${\bf k}=(k_1,\ldots,k_m)\in\mathbb{N}^m$ with $0\leq k_i\leq n_i$;

\item
Let a subword of $w_1\cdots w_{\beta-\alpha+|{\bf n}|}$ be a permutation of
$\{b_1^{n_1-k_1},\ldots,b_m^{n_m-k_m}\}$, with each $b_i$ weighted $y_i$ and also circled;
\item
Let all the other $w_i$'s be a permutation of $\{a^{\beta},b_1^{k_1},\ldots,b_m^{k_m}\}$, with each $b_i$ weighted
$x_i$ or $y_i$ and each $a$ weighted $1$.
\end{list}

We call such a weighted words $w$ a {\it configuration},
and define its weight as the product of the weights of all the $w_i$'s. It is not hard to see that the sum of the weights of the
configurations defined above is equal to ${\beta-\alpha+|{\bf n}| \choose {\bf n}-{\bf k}}{\beta+|{\bf k}| \choose
{\bf k}}({\bf x}+{\bf y})^{{\bf k}}{\bf y}^{{\bf n}-{\bf k}}$ for any ${\bf k}\in\mathbb{N}^m$.

Let $S^+$ and $S^-$ denote the sets of configurations with an even or
odd number of circled positions, respectively. We pair members $w$ of
$S^+$ and $S^-$ by identifying the first $w_j$ with weight $y_i$ ($1\leq i\leq m$) in $w_1\cdots w_{\beta-\alpha+|{\bf n}|}$,
either circling it or
uncircling it. Here we give an example for
 $\beta=4,\alpha=2$ and ${\bf n}=(2,2)$:

\begin{center}
\begin{picture}(360,20)(0,-10)
%\linethickness{1pt}
\put(0,0){$a$}\put(20,0){$b_2$} \put(42,0){$a$} \put(62,0){$b_1$}
\put(85,0){$b_2$} \put(108,0){$a$}
\put(128,0){$b_1$}\put(150,0){$a$}
\put(20,-15){$x_2$}\put(62,-15){$y_1$}\put(85,-15){$y_2$}\put(89,3){\circle{16}}
\put(128,-15){$x_1$} \put(164,0){$\Longleftrightarrow$}
\put(195,0){$a$}\put(215,0){$b_2$}
\put(215,-15){$x_2$}\put(237,0){$a$}
\put(257,0){$b_1$}\put(257,-15){$y_1$}
\put(280,0){$b_2$}\put(280,-15){$y_2$} \put(303,0){$a$}
\put(325,0){$b_1$}\put(325,-15){$x_1$}
\put(347,0){$a$}\put(261,3){\circle{16}}\put(284,3){\circle{16}}
\end{picture}
\end{center}

This pairs all configurations except those where every position
weighted $y_i$ is in the right $\alpha$ positions. These
configurations $w$ must belong to $S^+$ and the sum of their weights is equal to the right-hand side of \eqref{eq3}.
This is because if the subwords with elements weighted $y_i$ ($1\leq i\leq m$)
in $w_{\beta-\alpha+|{\bf n}|+1}\cdots w_{\beta+|{\bf n}|}$
is a permutation of $\{b_1^{n_1-k_1},\ldots,b_m^{n_m-k_m}\}$,
then there are ${\beta+|{\bf k}|\choose {\bf k}}$ possible ways
to choose the remaining subwords of $w$, where each $b_i$ is weighted $x_i$.
This proves \eqref{eq3}.


\begin{thebibliography}{99}
\small \setlength{\itemsep}{-.8mm}

\bibitem{Chapman}
R. Chapman,~A curious identity revised, Math. Gazette 87 (2003),
139--141.
\bibitem{Chen}
W. Y. C. Chen and S. X. M. Pang,~On the combinatorics of the Pfaff
identity, Discrete Math. 309 (2009) 2190--2196.
\bibitem{Guo}
V. J. W. Guo, Simple proofs of Jensen's, Chu's, Mohanty-Handa's,
and Graham-Knuth-Patashnik's identities, preprint, 2010, arXiv:1005.2745v1.
\bibitem{Munarini}
E. Munarini,~Generalization of a binomial identity of Simons,
Integers 5 (2005), \# A15.
\bibitem{Prodinger}
H. Prodinger,~A curious identity proved by Cauchy's integral
formula, Math. Gazette 89 (2005), 266--267.
\bibitem{Shattuck}
M. Shattuck,~Combinatorial proofs of some Simons-type binomial
coefficient identities, Integers 7 (2007), \# A27.
\bibitem{Simons}
S. Simons,~A curious identity, Math. Gazette 85 (2001),296--298.
\bibitem{Sun}
X. Wang and Y. Sun, A new proof of a curious identity, Math. Gazette
91 (2007), 105--106.
\end{thebibliography}
\end{document}